\documentclass[12pt]{amsart}
\usepackage[T1]{fontenc}
\usepackage{amssymb}
\newcommand{\Hee }{\mathcal H}

\newcommand{\Qee }{\mathcal Q}

\newcommand{\Pee }{\mathcal P}
\newcommand{\Fee }{\mathcal F}

\newtheorem{theorem}{Theorem}
\newtheorem{corollary}[theorem]{Corollary}
\newtheorem{lemma}[theorem]{Lemma}

\author{Piotr Kalemba}
\address{Piotr Kalemba \\
 Institute of Mathematics, University of
Silesia \\
ul. Bankowa 14, 40-007 Katowice}
\email{pkalemba@math.us.edu.pl}
\author{Szymon Plewik}
\address{Szymon Plewik\\Institute of Mathematics,
University of Silesia, ul. Ban\-ko\-wa 14, 40-007 Katowice}
\email{plewik@math.us.edu.pl}
\begin{document}

\title{Hausdorff gaps reconstructed from Luzin gaps}
\date{}
\begin{abstract}
  We consider a question: Can  a given AD-family be ADR for two orthogonal uncountable towers? If $b > \omega_1$, then we  rebuilt any AD-family of the cardinality $\omega_1$ onto a Hausdorff pre-gap.  Moreover, if a such AD-family is a Luzin gap, then we obtain a  Hausdorff gap. Under $b = \omega_1$,  a similar rebuilding is impossible.  
\end{abstract}

\maketitle
\section{Introduction}
 A family $\Qee$ is called \textit{almost disjoint}, briefly AD-\textit{family}, whenever any two members of $\Qee$ are almost disjoint, i.e. their intersection is finite. A set $C$ \textit{separates} a family $\Qee$ from a family $\Hee$, whenever each member  of $\Qee$ is almost contained  in $C$, i.e.  $B\setminus C$ is finite for any $B \in \Qee$,  and each member of $\Hee$ is almost disjoint with $C$. Whenever  sets $A$ and $B$ are almost disjoint for any $A\in \Qee$ and $B \in \Hee $, then families  $\Qee$ and $\Hee $ are called \textit{orthogonal}.  If  no set $C$ separates $\Qee$ from $\Hee$, then families  $\Qee$ and $\Hee $ are called \textit{non-separated}. Below, $A \subset^* B$ means that $A$ is almost contained in $B$, but not conversely.
   A pair of indexed families $[\{ A_\alpha: \alpha < \omega_1\} ; \{ B_\alpha: \alpha < \omega_1\}]$ is called \textit{Hausdorff pre-gap}, whenever 
$\alpha < \beta <\omega_1$ implies $ A_\alpha \subset^* A_\beta \subset^* B_\beta \subset^* B_\alpha $.  A Hausdorff pre-gap $[\{ A_\alpha: \alpha < \omega_1\} ; \{ B_\alpha: \alpha < \omega_1\}]$ is called \textit{Hausdorff gap}, whenever orthogonal towers $\{ A_\alpha: \alpha < \omega_1\}$ and $\{ \omega \setminus B_\alpha: \alpha < \omega_1\}$ are non-separated.  Establish that, a  family $\{ A_\alpha: \alpha < \lambda\}$ is a  \textit{tower}, whenever $\alpha < \beta$ implies $A_\alpha \subset^* A_\beta$.  
An AD-family $\Qee$ of the cardinality $\omega_1$ is called  \textit{Luzin gap}, whenever no two disjoint uncountable subfamilies of $\Qee$ are separated.  
An AD-family$ \Qee $  is  \textit{almost disjoint refinement} of a family $\Pee$ (briefly  $\Qee$ is ADR of $\Pee$), whenever   there exists a bijection $f: \Qee \to \Pee$ such that $X$ is almost contained in $f(X)$  for every
$X \in \Qee$.  Our definition of ADR is equivalent to the one considered in \cite{sim}, where  one can find a comprehensive discussion about almost disjoint refinements.

  We are going to compare constructions of Hausdorff and Luzin gaps. If $b > \omega_1$, then we describe how one can  rebuilt a AD-family of the cardinality $\omega_1$ onto a Hausdorff pre-gap. If a such AD-family is a Luzin gap, then we  obtain a  Hausdorff gap. Under $b = \omega_1$,  a similar rebuilding is impossible. For the sake of completeness, we enclose a construction of a Hausdorff gap which use no  form of so called the second interpolation theorem, compare \cite{sch}, and needs the hypothesis $b=\omega_1$.

 P. Simon indicated to us that Hausdorff gaps and Luzin gaps do not look compatible, September 2008 in Katowice.  M. Scheepers discerned something similar in \cite{sch}. Albeit,  he wrote that Luzin gaps are reminiscent of Hausdorff gaps. In \cite{kun}, K. Kunen declared that  "The easiest to construct are Luzin gaps" and that constructions of  Hausdorff gaps need some stronger inductive hypotheses.  Constructions of  Hausdorff gaps and Luzin gaps are considered apart, usually.    Hausdorff gaps have been examined via topological manner, through gap spaces associated with them, for example \cite{bs}, \cite{dou}  or \cite{nv}.  Forcing methods yield other treads to examine variety of Hausdorff gaps, for example \cite{as}, \cite{hir}, \cite{far}, \cite{sch} or \cite{yor}.  
  
\section{AD-families of the cardinality $b$.}    
 
 Recall that, $b$ is the least cardinality of  unbounded families of functions  $f: \omega \to \omega$ with respect to the partial order $\leq^*$, where $f\leq^* g$ whenever $f(n) \leq g(n)$ for all but finitely many $n\in \omega$. 
A function $h$ \textit{dominates} a restriction $f|_D$, whenever  $f(n) \leq h(n)$ for all but finitely many $n\in D$. If  $D=\omega$, then $h$ dominates $f$. 
It is well known that each of hypotheses $b=\omega_1$ or $b > \omega_1$ is consistent with ZFC. The hypothesis $b>\omega_1$ is  equivalent with  Proposition (1): \textit{The family of all sets of n.n. does not contain any $(\Omega, \omega^*)$ gaps}; by  Rothberger \cite{rot}. 
Consider the following question.  

\noindent \textbf{Question.}
 \textit{Could a given almost disjoint family be an almost disjoint refinement for the union of some two uncountable and orthogonal towers? } 
 
 \noindent To answer the question, we start with a  ZFC result. Then a Rothberger lemma is adapted in order to conclude some consistent results.  
 
\begin{theorem}\label{CC} There exists an almost disjoint family of the cardinality $b$, which is not almost disjoint refinement for any union of two orthogonal towers, where both towers have the  cardinality $b$.
\end{theorem}   

\begin{proof}

Let   $\Qee = \Fee \cup \{B_n: n < \omega\}$ be an AD-family such that always  $B_n= \{(n,k): k < \omega\}$ and  $\Fee = \{f_\alpha: \alpha < b\}$ consists of  almost disjoint and increasing functions   $f_\alpha: \omega \to \omega$. Assume that, $\Fee$ is unbounded and increasing. So, $\Qee$ consists of subsets of $ \omega \times \omega$ and  every $\Hee \subseteq \Fee$ of the cardinality $b$ is an unbounded family with respect to $\leq^*$.
  
 Suppose that  $\Qee$ is ADR of   the union of orthogonal towers $\{ A_\alpha: \alpha < b\}$ and $\{C_\alpha: \alpha < b\}$. Without loss of generality, one can fix $\alpha$ such that $C_\alpha$ almost contains infinitely many $B_n$.  Thus the family $$\Hee = \{f_\beta \in \Fee: f_\beta \subset^* \omega \times \omega \setminus C_\alpha\} $$ contains a subfamily $\Pee$ of the cardinality $b$ such that  $\Pee$ is an ADR of some subfamily of $\{ A_\alpha: \alpha < b\}$. So, the family  $\Hee$ is unbounded.
   On the other hand, put $h(n) = \max \{k: (n,k)\notin C_\alpha\}$ whenever  $B_n \subset^*C_\alpha$. Thus the function $h$ dominates each restriction $f_{\beta}|D$, where $f_\beta\in\Hee$ and $D=\{n: B_n\subset^* C_\alpha\}$.  Let $k_0, k_1, \ldots $ be an increasing enumeration of all elements of $D$. Put $g(i) =h(k_n)$ whenever   $k_{n-1}<i\leq k_n$. Because of $\Hee$ consists of increasing functions, one can check that $g$ dominates any function from $\Hee$; a contradiction. 
 \end{proof}

The following lemma can be derived from Rothberger's Lemma 5 stated in \cite{rot}.  
  
 \begin{lemma}\label{AA} Suppose  a countable family $\Qee$ consists of almost disjoint infinite subsets of natural numbers, and let $\Hee$ consists of sets almost disjoint with members of $\Qee$. If  $|\Hee|< b$, then families  $\Qee$ and $\Hee $ are separated. 
\end{lemma} 
\begin{proof} Without loss of generality, assume that members of $\Hee$ and $\Qee$ are subsets of $\omega \times \omega$ such that  $$\Qee =\{\{ (n,i): i \in \omega \}: n \in \omega \}.$$  Put $ f_B(n)=  \max \{i: (n,i)\in B \}$  for each $B\in \Hee $ (here $\max \emptyset =0$). Functions $f_B:\omega \to \omega$ are well defined since members of $\Hee$ are almost disjoint with elements of $\Qee$.  The family of all functions  $f_B$ has the cardinality less than $b$, so there exits a function $h$ which dominates each $f_B$. The set $$ \{(n,i): i>h(n) \mbox{ and } n\in \omega\}$$ separates $\Qee$ from $\Hee$.
\end{proof}

 Below, $A\subseteq^* B$ means that $A$ is almost contained in $B$.

\begin{theorem} \label{BB} 
Assume that $b>\omega_1$. If $\{ E_\alpha: \alpha < \omega_1\} \cup \{ F_\alpha: \alpha < \omega_1\}$ is an   AD-family, then there exists a Hausdorff pre-gap $$[\{ A_\alpha: \alpha < \omega_1\} ; \{ B_\alpha: \alpha < \omega_1\}]$$ such that $E_\alpha\subseteq^*A_{\alpha+1}\setminus A_\alpha\subseteq^*E_\alpha$ and $F_\alpha\subseteq^*B_{\alpha}\setminus B_{\alpha+1}\subseteq^*F_\alpha$, whenever $\alpha < \omega_1$.   
\end{theorem}
\begin{proof} We shall construct a desired  Hausdorff pre-gap, defining by induction sets $A_\alpha$ and $B_\alpha$ such that
\begin{enumerate}
	\item If $\beta < \alpha$, then $A_\beta \subset^* A_\alpha \subset^* B_\alpha \subset^* B_\beta$;
	\item If $\alpha = \beta +1$, then $ E_\beta \cup A_\beta = A_{\alpha}$  and   $    B_{\alpha} = B_\beta \setminus F_\beta   $;
	\item  Each member of the union  $\{ E_\beta: \alpha \leq \beta\} \cup \{ F_\beta: \alpha \leq \beta\}$ is almost disjoint with $A_\alpha$;
	\item  Each member of    $\{ E_\beta: \alpha \leq \beta\} \cup \{ F_\beta: \alpha \leq \beta\}$ is almost contained in $B_\alpha$.
\end{enumerate}

Put $A_0=\emptyset$ and $B_0=\omega$ and $A_{\alpha+1}= E_\alpha \cup A_\alpha $ and $B_{\alpha+1} = B_\alpha \setminus F_\alpha  $. It remains to  define sets $A_\alpha$ and $B_\alpha$
for limit ordinals $\alpha$.  Take  a sequence of ordinals $\gamma_0,  \gamma_1, \ldots$ which is  increasing and  has   the limit $\alpha$. Assume that  $\gamma_0=0$. 
 
At the first step,  let $ \Qee = \{A_{\gamma_{n +1}}\setminus A_{\gamma_n}: n \in \omega\}  $ and $ \Hee  = \{B_{\gamma_n}\setminus B_{\gamma_{n +1}}: n \in \omega\} \cup \{ E_\beta: \alpha \leq \beta\} \cup \{ F_\beta: \alpha \leq \beta\} .$ Families $\Qee$ and $\Hee$ are orthogonal and $\Qee$ is a countable AD-family. By  Lemma \ref{AA},   let $A_\alpha$ be a set which separates $\Qee$ from $\Hee $. Observe that  $\beta < \alpha$ implies $A_\beta \subset^* A_\alpha  \subset^* B_\beta$. Indeed, $\emptyset= A_{\gamma_{0}} \subset^* A_\alpha  \subset^* B_{\gamma_{0}}=\omega$. Inductively,  $A_{\gamma_{n}}\subseteq^*  (A_{\gamma_{n}}\setminus A_{\gamma_{n-1}}) \cup A_{\gamma_{n-1}}\subset^* A_\alpha $, since $A_\alpha$ separates $\Qee$ from $\Hee $. There exists $\gamma_{n} > \beta $,  hence $A_\beta \subset^* A_{\gamma_{n}} \subset^* A_\alpha$.  Also, one can assume that $A_\alpha \subset^* B_{\gamma_m}$. But sets $A_\alpha$ and $B_{\gamma_m}\setminus B_{\gamma_{m+1}}$ are almost disjoint, hence $A_\alpha \subset^* B_{\gamma_{m+1}}$. This gives that $A_\alpha \subset^* B_{\beta}$. 

At the second step, apply Lemma \ref{AA} to families  $\Qee = \{B_{\gamma_{n}}\setminus B_{\gamma_{n+1}}: n \in \omega\}$ and 
 $\Hee  = \{A_\alpha\} \cup \{ E_\beta: \alpha \leq \beta\} \cup \{ F_\beta: \alpha \leq \beta\} .$  Let $B_\alpha$ be the complement of a  set which separates $\Qee$ from $\Hee $, i.e. $B_\alpha$ separates $\Hee$ from $\Qee $. The union $\{B_\alpha\} \cup \{B_{\gamma_{n}}\setminus B_{\gamma_{n+1}}: n \in \omega\}$ is  an AD-family, hence   $\beta < \alpha $ implies $B_\alpha \subset^* B_\beta$.   
 \end{proof}
 
Thus,   one can reconstruct a Hausdorff gap from a Luzin gap, under  $b>\omega_1$. Indeed, let $\{ E_\alpha: \alpha < \omega_1\}$ and $ \{ F_\alpha: \alpha < \omega_1\}$ be    AD-families which are orthogonal and not separated. Then any  Hausdorff pre-gap like 
 in the Theorem \ref{BB}, i.e.   $[\{ A_\alpha: \alpha < \omega_1\} ; \{ B_\alpha: \alpha < \omega_1\}]$ such that $E_\alpha\subseteq^*A_{\alpha+1}\setminus A_\alpha\subseteq^*E_\alpha$ and $F_\alpha\subseteq^*B_{\alpha}\setminus B_{\alpha+1}\subseteq^*F_\alpha$, has to be a Hausdorff gap.   
 If  we assume that $\{ E_\alpha: \alpha < \omega_1\} \cup \{ F_\alpha: \alpha < \omega_1\}$ is a Luzin gap, then we have a  construction  of a  Hausdorff gap with some additional properties. 
 
 Let us  recall Luzin's construction of a gap, see \cite{luz}.  To convince the readers of Kunen's opinion, which is quoted in Introduction, we run as follows. Start with a family  $\{A_n:n\in \omega\}$ which consists of disjoint and infinite subsets of $\omega$.   Assume that almost disjoint sets $\{A_\beta: \beta < \alpha\}$  are just defined for a countable ordinal number $\alpha < \omega_1$. Enumerate  these sets $A_\beta$  into a sequence $ \{ B_n:n\in \omega\}$.  For every $n$, choose  a set  $$\{d_1, d_2, \ldots d_n\} \subset B_n\setminus (B_0\cup B_1 \cup \ldots \cup B_{n-1}),$$ with exactly  $n$ elements.  Than, put $A_\alpha$ to be the union of all already chosen sets $\{d_1, d_2, \ldots d_n\}$. The family $\{A_\alpha: \alpha < \omega_1\}$ forms  a Luzin gap. Indeed, consider  a partition of $\{A_\alpha: \alpha < \omega_1\}$ into two uncountably subfamilies $\mathcal{D}$ and $\mathcal{E}$. Suppose that  a set  set $C$ separates $\mathcal{D}$ from $\mathcal{E}$.  Fix  a natural number $n$ and uncountable subfamilies $\mathcal{F}\subseteq \mathcal{D}$ and $\mathcal{H}\subseteq \mathcal{E}$  such that $\cup \mathcal{F}\setminus n \subseteq C$ and $\cup \mathcal{H} \cap C \subseteq n$. Take $\alpha < \omega_1$ such that the intersection $\{A_\beta: \beta < \alpha \} \cap \mathcal{H}$ is infinite. Finally,  for each $\gamma > \alpha$ with $A_\gamma \in \mathcal{F}$  there exist $\beta<\alpha$ and  $A_\beta\in \mathcal{H}$ such that the intersection  $A_\beta \cap A_\gamma$ is a set $\{d_1, d_2, \ldots d_m\}$, where $m>n$.   This   is in conflict with   $\cup \mathcal{F}\setminus n \subseteq C$ and $\cup \mathcal{H} \cap C \subseteq n$.

 If  $b>\omega_1$ and there exists a Lebesgue non-measurable set of the cardinality $\omega_1$, then there exist AD-families of the cardinality $\omega_1$ which are non-measurable sets with respect to some Borel measures on $[\omega]^\omega$. But, any  family of sets which consists of a Hausdorff gap  has to be universally  measure zero, see \cite{ple}.  Thus, Hausdorff gaps and Luzin gaps could have consistently different measurable properties.

  \section{On constructions of Hausdorff gaps under $b=\omega_1$ } 
  
  It is consistent that   any AD-family of the cardinality $\omega_1$ is ADR of the union of some two orthogonal towers of the cardinality $\omega_1$ because of Theorem \ref{BB}. It is also clear  that this statement implies  $b>\omega_1$, since Theorem \ref{CC}  points out a suitable  AD-family. So,  we obtain a characterization of the hypothesis $b=\omega_1$.
 
\begin{corollary}  $b=\omega_1$  is equivalent with the existence of AD-family of the cardinality $\omega_1$ which is not an ADR of the union of any two orthogonal towers each of the cardinality $\omega_1$. \hfill $\Box$
 \end{corollary}

All known to us constructions of a Hausdorff gap use some forms of so called \textsl{The second interpolation theorem}, compare \cite{bs}, \cite{hau}, \cite{sie},\cite{sch} or \cite{yor}. In the previous part we do not use  this principle in inductive hypotheses.  So, we should add  constructions which  use no form of the second interpolation theorem. We use the following abbreviations: $\Delta =\{ (n,k)\in \omega \times \omega : k<n\}$ and $\int{f} =\{(n,k)\in \omega \times \omega: k\leq f(n)\}$.

Assume that  $b=\omega_1$. Let  $\{ \omega \setminus T_\alpha: \alpha < \omega_1 \}$ be a maximal tower and let  $\{ f_\alpha: \alpha < \omega_1\}$ be a unbounded family of functions, where $f_\alpha :\omega \to \omega$.
Let $A_0=\emptyset$ and $B_0= \omega \times \omega$ and fix a function $g_0:\omega \to \omega$ such that $ f_0  \leq^* g_0$.  Suppose that sets $A_\beta$ and $B_\beta$ and functions $g_\beta$ are defined for $\beta < \alpha$. We should define sets $A_\alpha$ and $B_\alpha$ and a function $g_\alpha$ such that 
\begin{itemize}
	\item  $ A_\beta \subset^* A_\alpha \subset B_\alpha \subset^* B_\beta $ for each $\beta < \alpha$;
	\item $g_\alpha\subseteq A_\alpha\subseteq \int{g_\alpha}$, where a function $g_\alpha$ is  such that $ f_\alpha  \leq^* g_\alpha$;
	 \item $B_\alpha = A_\alpha \cup (\omega\times T_\alpha)\setminus \Delta$.
\end{itemize}
To do this, take $X$ such that $ A_\beta \subset^* X \subset^* B_\beta $ for each $\beta < \alpha$. Fix a function $g_\alpha \subset\omega\times T_\alpha \setminus \Delta$ such that $g_\alpha$ dominates every function from $\{g_\beta: \beta < \alpha\} \cup \{f_\alpha\}$. Eventually, put $A_\alpha = g_\alpha \cup (X\cap \int{g_\alpha})   $ and  $B_\alpha = A_\alpha \cup (\omega\times T_\alpha)\setminus \Delta$. Above defined sets $A_\alpha$ and $B_\alpha$ constitute a Hausdorff pre-gap. The tower    $\{ \omega \setminus T_\alpha: \alpha < \omega_1 \}$ is  maximal. Hence, whenever $A_\alpha \subset^* C\subset^*  B_\alpha$ for any $\alpha < \omega_1$, then there exists a function $h$ such that $C \subset \int{h}$. But this means that $h$ dominates each $f_\alpha$, a contradiction.


\begin{thebibliography}{40}
\bibitem{as} U. Abraham, S. Shelah, \textit{Ladder gaps over stationary sets}, J. Symbolic Logic 69 (2004),  518 - 532.
\bibitem{bs} A. B\l aszczyk, A. Szyma\'nski, 
\textit{Hausdorff's gaps versus normality} 
Bull. Acad. Polon. Sci. Sér. Sci. Math. 30  no. 7-8,(1982), 371--378. 
\bibitem{dou} E. van Douwen, 
\textit{Hausdorff gaps and a nice countably paracompact nonnormal space},  Topol. Proc., Vol. 1 (1976), 239 - 242. 
\bibitem{far} I. Farah, 
\textit{Luzin gaps},  
Trans. Am. Math. Soc. 356, No. 6 (2004) 2197 - 2239.
\bibitem{hau} F. Hausdorff, \textit{Summen von $\aleph_1$ Mengen}, Fund. Math. 26 (1936), 243 - 247. 
\bibitem{hir} J. Hirschorn, \textit{On the strenght of Hausdorff's gap condition}, arXiv:0806.4732v1.
\bibitem{luz} N. Luzin, \textit{On subsets of the series of natural numbers}, Isv. Akad. Nauk. SSSR Ser. Mat. 11
(1947), 403 - 411. 
\bibitem{kun} K. Kunen, \textit{Where MA first fails},  
J. Symbolic Logic 53 (1988), no. 2, 429 - 433.
\bibitem{nv} P. J. Nyikos and J. E. Vaughan,  \textit{On first countable, countably compact spaces. I. $(\omega_1 ,\omega_1^*)$-gaps}, Trans. Amer. Math. Soc. 279 (1983), no. 2, 463 - 469.
\bibitem{ple} Sz. Plewik, 
\textit{Towers are universally measure zero and always of first category}, 
Proc. Amer. Math. Soc. 119 (1993), no. 3, 865 - 868. 
\bibitem{rot} F. Rothberger,  \textit{On some problems of Hausdorff and of Sierpi\'nski}, Fund. Math. 35, (1948), 29 - 46.
\bibitem{sch} M. Scheepers, 
\textit{Gaps in $(\sp\omega  \omega, \prec)$},  Set theory of the reals (Ramat Gan, 1991), 439--561, 
Israel Math. Conf. Proc., 6, Bar-Ilan Univ., Ramat Gan, (1993).
\bibitem{sie} W. Sierpi\'nski,  \textit{General topology}, Translated by C. Cecilia Krieger. Mathematical Expositions, No. 7, University of Toronto Press, Toronto, (1952).
\bibitem{sim} P. Simon,  \textit{A note on almost disjoint refinement}, 24th Winter School on Abstract Analysis (Bene$\check{\mbox{s}}$ova Hora, 1996). Acta Univ. Carolin. Math. Phys. 37 (1996), no. 2, 89 - 99. 
\bibitem{yor} T. Yorioka, \textit{Some results on gaps in $\Pee(\omega)/$fin}, unpublished?
\end{thebibliography}
\end{document}